\documentclass[12pt,a4paper]{article}
\usepackage {amssymb,latexsym}
\usepackage {amsmath}

\newtheorem{theorem}{Theorem}[section]
\newtheorem{proposition}{Proposition}[section]
\newtheorem{remark}{Remark}[section]
\newtheorem{proof}{Proof}

\newtheorem{definition}{Definition}

\newcommand{\real}{\mathbb R}
\newcommand{\corps}{\mathbb K}

\newcommand{\complex}{\mathbb C}

\newcommand{\B}{\mathbb B}

\newcommand{\ben}{\begin{enumerate}}
\newcommand{\een}{\end{enumerate}}

\newcommand{\A}{\mathcal{A}}

\begin{document}
\author{Abdenacer MAKHLOUF} 

\title{A comparison of deformations and geometric study of associative algebras
varieties}

\maketitle
\begin{center}
Laboratoire de Math\'{e}matiques et Applications, \\
Universit\'{e} de Haute Alsace, Mulhouse, France \vspace{4mm}
\\Email :\\
 Abdenacer.Makhlouf@uha.fr 
\end{center}

{\it The aim of this paper is to give an overview and to compare the different deformation
theories of algebraic structures. We describe in
each case the corresponding notions of degeneration and rigidity. We
illustrate these notions with examples and give some general properties. The last part of this work shows how these notions help in the
study of associative algebras varieties.
The first and popular deformation approach was introduced by M. Gerstenhaber for rings and 
algebras using formal power series. 
A noncommutative version was given by Pinczon and
generalized by F. Nadaud. A more general approach called global deformation
follows from a general theory of Schlessinger and was developed by A. Fialowski in order to deform 
infinite dimensional nilpotent Lie algebras. In a Nonstandard
framework, M.Goze introduced the notion of perturbation for studying the rigidity
of finite dimensional complex Lie algebras. 
All these approaches has in common the fact that we make an
"extension" of the field. These theories may be applied to any multilinear structure. We shall be concerned in this paper with the category of associative
algebras.}

AMS Classification 2000 : 16-xx, 13-xx, 14-Rxx, 14-Fxx, 81-xx.

\section{Introduction}

Throughout this paper $\corps$ will be an algebraically closed field, and $\mathcal{A}$
denotes an associative $\corps$-algebra. Most examples will be in the finite 
dimensional case, let $V$ be the underlying 
$n$-dimensional vector space of $\mathcal{A}$ over $\corps$ and $(e_{1},\cdots ,e_{n})$ be a basis
of $V$. The bilinear map $\mu $ denotes the multiplication of $\mathcal{A}$ on $V$,
 and  $e_{1}$ is 
the unity. By linearity this can be done by specifying the $n^{3}$
structure constants $C_{ij}^{k}\in \corps$ where $\mu
(e_{i},e_{j})=\sum_{k=1}^{n}C_{ij}^{k}e_{k}$. The associativity condition
 limits the sets of structure constants, 
 $C_{ij}^{k}$ to a subvariety of $\corps ^{n^3}$ which we denote by $Alg_{n}$. 
 It is generated by the polynomial relations
 \[
\sum_{l=1}^{n}C_{ij}^{l}C_{lk}^{s}-C_{il}^{s}C_{jk}^{l}=0\ and\
C_{1i}^{j}=C_{i1}^{j}=\delta _{ij}\quad 1\leq i,j,k,s\leq n. 
\]
This variety is quadratic, non regular and in general nonreduced. The
natural action of the group $GL(n,\corps)$ corresponds to the change of basis :
two algebras  $\mu _{1\text{ }}$and $\mu _{2}$ over $V$ are isomorphic if there exists 
$f$ in $GL(n,\corps)$ such that : 
\[
\forall X,Y\in V\quad \mu _{2}\left( X,Y\right) =\left( f\cdot \mu
_{1}\right) \left( X,Y\right) =f^{-1}\left( \mu _{1}\left( f\left( X\right)
,f\left( Y\right) \right) \right) 
\]
The orbit of an algebra $\A$ with multiplication $\mu _{0}$, denoted by $\vartheta \left( \mu
_{0}\right) $, is the set of all its isomorphic algebras. The deformation techniques are used to
 do the geometric study of these varieties. The deformation attempts to understand which algebra 
 we can get from the original
 one by deforming. At the same time it gives more  informations 
 about the structure of the algebra, for example we can look for which properties are stable under
 deformation. 
   
The deformation of mathematical objects is one of the oldest techniques used
by mathematicians.  The different areas where the notion of deformation
appears are geometry, complex manifolds (Kodaira and Spencer 1958, Kuranishi
1962), algebraic manifolds (Artin, Schlessinger 1968), Lie algebras
(Nijenhuis, Richardson 1967) and rings and associative algebras
(Gerstenhaber 1964). The dual notion of deformation (in some sense) is that
of degeneration which appears first in physics literature ( Segal 1951, Inonu and Wigner
1953). Degeneration is also called specialisation or contraction.

The quantum mechanics and quantum groups gave impulse to the theory of deformation. The 
theory of quantization was introduced by Bayen, Flato, Fronsdal, Lichnerowicz and Sternheimer \cite{BFLS} to describe
the quantum mechanics as a deformation of classical mechanics. It associates to a Poisson manifold a star-product, which is a one
parameter family  of associative algebras. In another hand the quantum groups are
obtained by deforming the Hopf structure of an algebra, in particular enveloping algebras. In \cite{BMP},
we show the existence of associative deformation of an enveloping algebra using the linear Poisson
structure of the Lie algebra.

In the following, we will recall first the different notions of deformation.
The most used one is the formal deformation introduced by Gerstenhaber for
rings and algebras \cite{Gerstenhaber64}, it uses a formal series and connects the
theory of deformation to Hochschild cohomology. A noncommutative 
version, where the parameter no longer commutes with the element of algebra, was introduced by Pinzson \cite{pincson} and generalized
by Nadaud \cite{Nadaud1}. They describe the corresponding cohomology and
show that the Weyl algebra, which is rigid for formal deformation, is nonrigid 
in the noncommutative case. In a Nonstandard framework, Goze
introduced the notion of perturbation for studying the rigidity of Lie
algebras\cite{Goze88}\cite{AncocheaGoze02}, it was also used to describe the
6-dimensional rigid associative algebras \cite{Goze-Makhlouf90}\cite
{Makhlouf-Goze96}. The perturbation needs the concept of infinitely small elements, these
elements are obtained here in an algebraic way. We construct an extension of real or
complex numbers sets containing these elements.
  A more general notion called global deformation was
introduced by Fialowski, following Schlessinger, for Lie algebras \cite
{Fialowski88}. All these approaches are compared in section 3. The section 4 is devoted to
 universal and versal deformation. In Section 5, we describe the corresponding notions of
degeneration with some general properties and examples. The Section 6 is
devoted to the study of rigidity of algebra in each framework. The last section
concerns the geometric study of the algebraic varieties $A\lg _{n}$ using
the deformation tools.

In formal deformation the properties are described
using the Hochschild cohomology groups. The global deformation seems to be the good framework to solve the
problem of universal or versal deformations, the deformations which generate the others. 
Whereas, the perturbation approch is more adapted to direct computation.

\section{The world of deformations}

\subsection{Gerstenhaber's formal deformation }

Let $\A$ be an associative algebra over a field $\corps$, $V$ be the underlying
vector space and $\mu _{0}$ the multiplication. \newline
Let $\corps\left[ \left[ t\right] \right] $ be the power series ring in one
variable $t$ and $V\left[ \left[ t\right] \right] $ be the extension of $V$
by extending the coefficient domain from $\corps$ to $\corps \left[ \left[ t\right]
\right] $. Then $V\left[ \left[ t\right] \right] $ is a $\corps\left[ \left[
t\right] \right] $-module and $V\left[ \left[ t\right] \right] =V\otimes
_{\corps}\corps\left[ \left[ t\right] \right] $. Note that $V$ is a submodule of $%
V\left[ \left[ t\right] \right] $. We can obtain an extension of $V$ with a
structure of vector space by extending the coefficient domain from $\corps$ to $%
\corps\left( \left( t\right) \right) $, the quotient power series field of $%
\corps\left[ \left[ t\right] \right] $.\newline
Any bilinear map $f:V\times V\rightarrow V$ (in particular the
multiplication in $\A$) can be extended to a bilinear map from $V\left[
\left[ t\right] \right] \times V\left[ \left[ t\right] \right] $ to $V\left[
\left[ t\right] \right] .$

\begin{definition}
{\it Let }$\mu _{0}$ {\it be the multiplication of the associative algebra }$%
\A$. {\it A deformation of }$\mu _{0}$ {\it is a one parameter family} $\mu
_{t}$ {\it in }$\corps\left[ \left[ t\right] \right] \otimes V${\it \ over the
formal power series ring }$\corps\left[ \left[ t\right] \right] $ {\it of the form%
} $\mu _{t}=\mu _{0}+t\mu _{1}+t^{2}\mu _{2}+....\qquad where\quad \mu
_{i}\in Hom\left( V\times V,V\right) $ (bilinear maps) satisfying the (formally) condition of associativity : 
\[
\text{For all }X,Y,Z\in V\quad \mu _t\left( \mu _t\left( X,Y\right)
,Z\right) =\mu _t\left( X,\mu _t\left( Y,Z\right) \right) 
\]
\end{definition}

We note that the deformation of $\A$ is a $\corps$-algebra structure on $\A [[t]]$ such that
$\A [[t]]/ t \A [[t]]$ is isomorphic to $\A$.

The previous equation is equivalent to an infinite equation system and it is 
 called the {\it %
deformation equation}. The resolution of the deformation equation connects deformation theory to Hochschild
cohomology.
Let $C^d ( \mathcal{A}, \mathcal{A} )$ be the space of $d$-cochains, the space of
 multilinear maps from $V^{\times d}$ to $V$. 
 
The boundary operator $\delta_{d}$, which we denote by $\delta$ if there is no ambiguity
:
$$
\delta_{d} : 
\begin{array}{c}
C ^d\left( \mathcal{A}, \mathcal{A}\right) \rightarrow C ^{d+1}\left( \mathcal{A}, \mathcal{A}\right) \\ 
\varphi \rightarrow \delta_{d}\varphi 
\end{array}
$$
is defined for $\left( x_1,...,x_{d+1}\right) \in V^{\times (d+1)}$ by

$
\begin{array}{c}
\delta_{d}\varphi \left( x_1,...,x_{d+1}\right) =\mu_{0} \left( x_1,\varphi
\left( x_2,...,x_{d+1}\right) \right) + \\ \sum_{i=1}^d\left( -1\right)
^i\varphi \left( x_1,...,\mu_{0} \left( x_i,x_{i+1}\right) ,...,x_d\right) + 
\left( -1\right) ^{d+1}\mu_{0} \left( \varphi \left( x_1,...,x_d\right)
,x_{d+1}\right) 
\end{array}
$

The  group of $d$-cocycles is :
$
Z^d\left( \A,\A\right) =\left\{ \varphi \in C ^d\left( \mathcal{A}, \mathcal{A}\right) /\
\delta ^d\varphi =0\right\} 
$

The group of $d$-coboundaries is 
\begin{center}
$
B^d\left( \A,\A\right) =\left\{ \varphi \in C ^d\left( \mathcal{A}, \mathcal{A}\right) /\
\varphi =\delta _{d-1}f,\ f\in C ^{d-1}\left( \A,\A\right) \right\} 
$
\end{center}
The Hochschild cohomology group of the algebra $\A$ with coefficient in
itself is

\begin{center}
$H^d\left( \A,\A\right) =Z^d\left( \A,\A\right) /B^d\left( \A,\A\right) .$
\end{center}

We define two maps $\circ$ and $[~,~]_{G}$
\begin{center}
$ \circ , [ ., . ]_G : C^d ( \mathcal{A}, \mathcal{A} )
   \times C^e ( \mathcal{A}, \mathcal{A} ) \rightarrow C^{d + e
   - 1} ( \mathcal{A}, \mathcal{A} ) $
\end{center}  
by 
\begin{center}
$( \varphi \circ \psi ) ( a_1, \cdots, a_{d + e - 1} ) = \sum_{i \geqslant
0} ( - 1 )^{i ( e - 1 )} \varphi ( a_1, \cdots, a_i, \psi ( a_{i + 1}, \cdots,
a_{i + e} ), \cdots )$ 

and 

$[ \varphi, \psi ]_G =\varphi \circ \psi - ( - 1
)^{( e - 1 ) ( d - 1 )} \psi \circ \varphi$.
\end{center} 
The space $(
C ( \mathcal{A}, \mathcal{A} ) ,\circ)$ is a  pre-Lie
algebra and $(
C ( \mathcal{A}, \mathcal{A} ) \text{$[ ., . ]_G$)}$ is a graded Lie
algebra (Lie superalgebra). The bracket $[ ., . ]_G$ is called Gerstenhaber's bracket. The square of $[
\mu_{0}, . ]_G$ vanishes and defines the $2$-coboundary operator. 
The multiplication $\mu_{0}$ of $\mathcal{A}$ is associative if $[ \mu_{0},
\mu_{0} ]_G = 0$.

Now, we discuss the deformation equation in terms of cohomology. The deformation
equation  may be written 
\[
\left( \ast \right) \quad \left\{ \sum_{i=0}^{k}\mu _{i}\circ \mu
_{k-i} =0 \qquad k=0,1,2,\cdots \right. 
\]
The first equation $\left( k=0\right) $ is the associativity condition for $%
\mu _{0}$. 
The second equation shows that $\mu _{1}$ must be a 2-cocycle for Hochschild
cohomology $\left( \mu _{1}\in Z^{2}\left( \A,\A\right) \right) $.\newline
More generally, suppose that $\mu _{p}$ be the first non-zero coefficient
after $\mu _{0}$ in the deformation $\mu _{t}$. This $\mu _{p}$ is called
the {\it infinitesimal} of $\mu _{t}$ and {\it is a} 2-{\it cocycle of the
Hochschild cohomology of }$\A$ {\it with coefficient in itself.}

{\it The cocycle }$\mu _{p}$ {\it is called integrable if it is the first
term, after }$\mu _{0},$ {\it of an associative deformation.}

The integrability of $\mu _p$ implies an infinite sequence of relations
which may be interpreted as the vanishing of the obstruction to the
integration of $\mu _p$.

For an arbitrary  $k>1,$ the $k^{th}$ equation of the system $\left(
*\right) $ may be written 
\[
\delta \mu _{k} =\sum_{i=1}^{k-1}\mu _{i}\circ \mu_{k-i}
\]
Suppose that the truncated deformation $\mu _{t}=\mu _{0}+t\mu _{1}+t^{2}\mu
_{2}+....t^{m-1}\mu _{m-1}$ satisfies the deformation equation. The
truncated deformation is extended to a deformation of order $m$, i.e. $\mu
_{t}=\mu _{0}+t\mu _{1}+t^{2}\mu _{2}+....t^{m-1}\mu _{m-1}+t^{m}\mu _{m}$, 
satisfying the deformation equation if 
\[
\delta \mu _{m}=\sum_{i=1}^{m-1}\mu _{i}\circ \mu_{m-i}
\]
The right-hand side of this equation is called the {\it obstruction }to
finding $\mu _{m}$ extending the deformation.

{\it The obstruction is a Hochschild 3-cocycle. Then, if }$H^{3}\left(
\A,\A\right) =0${\it \ it follows that all obstructions vanish and every }$\mu _{m}\in
Z^{2}\left( \A,\A\right) ${\it \ is integrable.}

Given two associative deformations $\mu _{t\text{ }}$and $\mu _{t}^{^{\prime
}}$ of $\mu _{0}$, we say that they are equivalent if there is a formal
isomorphism $F_{t}$ which is a $\corps\left[ \left[ t\right] \right] $-linear map
that may be written in the form 
\[
F_{t}=Id+tf_{1}+t^{2}f_{2}+\cdots \qquad where\quad f_{i}\in End_{\corps}\left(
V\right) 
\]
such that $\mu _{t}=F_{t}\cdot \mu _{t}^{^{\prime }}$ defined by 
\[
\mu _{t}(X,Y)=F_{t}^{-1}\left( \mu _{t}^{\prime }\left( (F_{t}\left( X\right)
,F_{t}\left( Y\right) \right) \right) \text{\quad for all }X,Y\in V 
\]
A deformation $\mu _{t\text{ }}$of $\mu _{0}$ is called  {\it trivial} if
and only if $\mu _{t\text{ }}$is equivalent to $\mu _{0}$.

\begin{proposition}
{\it Every nontrivial deformation } $\mu _{t}${\it \ of }$\mu _{0}${\it \ is equivalent
to }$\mu _{t}=\mu _{0}+t^{p}\mu _{p}^{\prime }+t^{p+1}\mu _{p+1}^{\prime
}+\cdots ${\it \ where }$\mu _{p}^{\prime }\in Z^{2}\left( \A,\A\right) ${\it %
\ and} $\mu _{p}^{\prime }\notin B^{2}\left( \A,\A\right) $.
\end{proposition}

Then we have this fundamental and well-known theorem.

\begin{theorem}
If $H^{2}\left( \A,\A \right) =0$, then all deformations of $\A$ are equivalent
to a trivial deformation.
\end{theorem}

\begin{remark}
The formal deformation notion is extended to coalgebras and bialgebras in \cite{Gerstenhaber90}.
\end{remark}

\subsection{Non-commutative formal deformation}

In previous formal deformation the  parameter commutes with the
original algebra. Motivated by some nonclassical deformation appearing in
quantization of Nambu mechanics, Pinczon introduced a deformation called
noncommutative deformation where the parameter no longer commutes with the
original algebra. He developed also the associated cohomology \cite{pincson}.

Let $\A$ be a $\corps$-vector space and $\sigma $ be an endomorphism of $\A$. We
give $\A[[t]]$ a $\corps [[t]]$-bimodule structure defined for every $a_{p}\in
\A,\lambda _{q}\in \corps$ by : 
\begin{eqnarray*}
\sum_{p\geq 0}a_{p}t^{p}\cdot \sum_{q\geq 0}\lambda _{q}t^{q}
&=&\sum_{p,q\geq 0}\lambda _{q}a_{p}t^{p+q} \\
\sum_{q\geq 0}\lambda _{q}t^{q}\cdot \sum_{p\geq 0}a_{p}t^{p}
&=&\sum_{p,q\geq 0}\lambda _{q}\sigma ^{q}(a_{p})t^{p+q}
\end{eqnarray*}

\begin{definition}
A $\sigma $-deformation of an algebra $\A$ is a $\corps$-algebra structure on $%
\A[[t]]$ which is compatible with the previous $\corps[[t]]$-bimodule structure
and such that $$\A[[t]]/(\A[[t]]t)\cong \A.$$
\end{definition}
The previous deformation was generalized by F. Nadaud \cite{Nadaud1} where
he considered deformations based on two commuting endomorphisms $\sigma $ and 
$\tau $. The $\corps[[t]]$-bimodule structure on $\A[[t]]$ is defined for $a\in \A$
by the formulas $t\cdot a=\sigma (a)t$ and $a\cdot t=\tau (a)t$, ($a\cdot t$
being the right action of $t$ on $a$).

The remarquable difference with commutative deformation is that the Weyl
algebra of differential operators with polynomial coefficients over $\real$ is
rigid for commutative deformations but has a nontrivial noncommutative
deformation; it is given by the enveloping algebra of the Lie superalgebra
$osp(1,2)$.

\subsection{Global deformation}
The approach follows from a general fact in Schlessinger's works \cite{schliss1} 
and was developed by A. Fialowski \cite{Fialowski86}. She applies it to construct
deformations of Witt Lie subalgebras. 
We summarize here the notion of global deformation in the case of associative algebra. Let $\B$ be a commutative 
algebra over a field $\corps$ of characteristic 0 and 
augmentation morphism $\varepsilon :\mathcal{A}\rightarrow \corps$ (a $\corps$-algebra homomorphism, $%
\varepsilon (1_{\B})=1$). We set $m_{\varepsilon }=Ker(\varepsilon )$; $%
m_{\varepsilon }$ is a maximal ideal of $\B$. (A maximal ideal $m$ of $%
\B $ such that $\A/m\cong \corps$, defines naturally an augmentation). We call $(\B,m)$
base of deformation.

\begin{definition}
A global deformation of base $(\B,m)$ of an algebra $\mathcal{A} $ with a
multiplication $\mu $ is a structure of $\B$-algebra on the tensor product $%
\B\otimes _{\corps}\mathcal{A} $ with a multiplication $\mu _{\B }$ such that $%
\varepsilon \otimes id:\B\otimes \mathcal{A} \rightarrow \corps\otimes \mathcal{A} =%
\mathcal{A} $ is an algebra homomorphism. i.e. $\forall a,b\in \B$ and $\forall
x,y\in \mathcal{A} $ :

\begin{enumerate}
\item  $\mu _{\B }(a\otimes x,b\otimes y)=(ab\otimes id)\mu _{\B
}(1\otimes x,1\otimes y)\quad (\B-$linearity$)$
\item  The multiplication $\mu _{\B }$ is associative.
\item  $\varepsilon \otimes id\left( \mu _{\B }(1\otimes x,1\otimes
y)\right) =1\otimes \mu (x,y)$
\end{enumerate}
\end{definition}

\begin{remark}
Condition 1 shows that to describe a global deformation it is enough to know
the products $\mu _{\B }(1\otimes x,1\otimes y)$, where $x,y\in \mathcal{A} .
$ The conditions 1 and 2 show that the algebra is associative and the last
condition insures the compatibility with the augmentation. We deduce 
\[
\mu _{B}(1\otimes x,1\otimes y)=1\otimes \mu (x,y)+\sum_{i}\alpha
_{i}\otimes z_{i}\quad \text{with }\alpha _{i}\in m,\ z_{i}\in \mathcal{A} 
\]
\end{remark}
\subsubsection{Equivalence and push-out}
\begin{itemize}
\item A global deformation is called {\em trivial\em} if the structure of $\B$-algebra on $%
\B\otimes _{\corps}\mathcal{A} $ satisfies $\mu _{\B }(1\otimes x,1\otimes
y)=1\otimes \mu (x,y).$
\item Two deformations of an algebra with the same base are called {\em equivalent \em} (or
isomorphic) if there exists an algebra isomorphism between the two copies of 
$\B\otimes _{\corps}\mathcal{A} $, compatible with $\varepsilon \otimes id.$
\item A global deformation with base $(\B,m)$ is called {\em local \em} if $\B$ is a local $\corps$%
-algebra with a unique maximal ideal $m_{\B}$. If, in addition $m_{\B}^{2}=0$, 
the deformation is called {\em infinitesimal \em}.
\item Let $\B^{\prime }$ be another commutative algebra over $\corps$ with augmentation $%
\varepsilon ^{\prime }:\B^{\prime }\rightarrow \corps$ and $\Phi :\B\rightarrow
\B^{\prime }$ an algebra homomorphism such that $\Phi (1_{\B})=1_{\B^{\prime }}$
and $\varepsilon ^{\prime }\circ \Phi =\varepsilon $. If a deformation $%
\mu_{\B} $ with a base $(\B,Ker(\varepsilon ))$ of $\mathcal{A} $ is given we call
push-out $\mu_{\B ^{\prime }}=\Phi _{*}\mu_{\B} $ a deformation of $%
\mathcal{A} $ with a base $(\B^{\prime },Ker(\varepsilon ^{\prime }))$ with the
following algebra structure on $\B^{\prime }\otimes \mathcal{A} %
=\left( \B^{\prime }\otimes _{\B}\B\right) \otimes \mathcal{A} =\B^{\prime }\otimes
_{\B}\left( \B\otimes \mathcal{A} \right) $ $$\mu _{\B ^{\prime }}\left( a_{1}^{\prime
}\otimes _{\B}\left( a_{1}\otimes x_{1}\right) ,a_{2}^{\prime }\otimes
_{\B}\left( a_{2}\otimes x_{2}\right) \right) :=a_{1}^{\prime }a_{2}^{\prime
}\otimes _{\B}\mu _{\B }\left( a_{1}\otimes x_{1},a_{2}\otimes
x_{2}\right) $$ with $a_{1}^{\prime },a_{2}^{\prime }\in \B^{\prime
},a_{1},a_{2}\in \B,x_{1},x_{2}\in \mathcal{A} $. The algebra $\B^{\prime }$ is viewed
as a $\B$-module with the structure $aa^{\prime }=a^{\prime }\Phi \left(
a\right)$. Suppose that  $$\mu _{\B }(1\otimes
x,1\otimes y)=1\otimes \mu (x,y)+\sum_{i}\alpha _{i}\otimes z_{i}\quad $$ with $\alpha _{i}\in
m,\ z_{i}\in \mathcal{A} $,  then $$\mu _{\B^{\prime
}}(1\otimes x,1\otimes y)=1\otimes \mu (x,y)+\sum_{i}\Phi (\alpha
_{i})\otimes z_{i}\quad $$.
\end{itemize}

\subsubsection{Coalgebra and Hopf algebra global deformation}
The global deformation may be extended to  coalgebra structures, then to Hopf algebras. 
Let $\mathcal{C}$ be a coalgebra over $\corps$, defined by the comultiplication $\Delta
:\mathcal{C}\rightarrow \mathcal{C}\otimes C$. Let $\B$ be a commutative algebra over $\corps$ and let
$\varepsilon$ be an
augmentation $\varepsilon :\B\rightarrow \corps$ with $m=Ker(\varepsilon )$ a
maximal ideal .

A global deformation with base $(\B,m)$ of coalgebra $C$ with a
comultiplication $\Delta $ is a structure of $\B$-coalgebra on the tensor
product $\B\otimes _{\corps}\mathcal{C}$ with the comultiplication $\Delta _{\B }$ such
that $\varepsilon \otimes id:\B\otimes \mathcal{C}\rightarrow \corps\otimes \mathcal{C}=\mathcal{C}$ is a
coalgebra homomorphism. i.e $\forall a\in \B$ and $\forall x\in \mathcal{C}$ :

\begin{enumerate}
\item  $\Delta _{\B }(a\otimes x)=a\Delta _{\B }(1\otimes x)$
\item  The comultiplication $\Delta _{\B }$ is coassociative.

\item  $\left( \varepsilon \otimes id\right) \otimes \left( \varepsilon
\otimes id\right) \left( \Delta _{\B }(1\otimes x)\right) =1\otimes
\Delta (x)$
\end{enumerate}

The comultiplication $\Delta _{\B }$ may be written for all $x$ 	as :
\[
\Delta _{B }(1\otimes x)=1\otimes \Delta (x)+\sum_{i}\alpha
_{i}\otimes z_{i}\otimes z_{i}^{\prime }\quad \text{with }\alpha _{i}\in m,\
z_{i},z_{i}^{\prime }\in \mathcal{C}
\]
\subsubsection{Equivalence and push-out for coalgebras}
Two global deformations of a coalgebra with the same base are called {\em equivalent \em} (or
isomorphic) if there exists a coalgebra isomorphism between the two copies of 
$\B\otimes _{\corps}\mathcal{C} $, compatible with $\varepsilon \otimes id.$

Let $\B^{\prime }$ be another commutative algebra over $\corps$ as in section 2.3.1. 
If  $\Delta_{\B} $ is a global deformation with a base $(\B,Ker(\varepsilon ))$ of 
a coalgebra $\mathcal{C} $. We call
push-out $\Delta_{\B ^{\prime }}=\Phi _{*}\Delta_{\B} $ a global deformation 
with a base $(\B^{\prime },Ker(\varepsilon ^{\prime }))$ of $\mathcal{C} $  with the
following coalgebra structure on $\B^{\prime }\otimes \mathcal{C} %
=\left( \B^{\prime }\otimes _{\B}\B\right) \otimes \mathcal{C} =\B^{\prime }\otimes
_{\B}\left( \B\otimes \mathcal{C} \right) $ $$\Delta_{\B ^{\prime }}\left( a^{\prime
}\otimes _{\B}\left( a \otimes x \right) \right) :=a^{\prime }\otimes _{\B}\Delta _{\B }\left( 
a\otimes x \right) $$ with $a^{\prime }\in \B^{\prime
},a \in \B,x \in \mathcal{C} $. The algebra $\B^{\prime }$ is viewed
as a $\B$-module with the structure $aa^{\prime }=a^{\prime }\Phi \left(
a\right)$. 
Suppose that  $$\Delta _{\B }(1\otimes
x)=1\otimes \Delta (x)+\sum_{i}\alpha _{i}\otimes z_{i}\otimes z^{\prime}_{i}\quad $$ with $\alpha _{i}\in
m,\ z_{i}\in \mathcal{C} $,  then $$\Delta _{\B^{\prime
}}(1\otimes x)=1\otimes \Delta (x)+\sum_{i}\Phi (\alpha
_{i})\otimes z_{i}\otimes z^{\prime}_{i}\quad $$.

\subsubsection{Hopf algebra global deformation}
Naturally, we can define Hopf algebra global deformation from the algebra and coalgebra global
deformation.

\subsubsection{Valued global deformation.}

In \cite{Goze-Remm02}, Goze and Remm considered the
case where the base algebra $\B$ is a commutative $\corps$-algebra of valuation
such that the residual field $\B /m$ is isomorphic to $\corps$, where $m$ is the
maximal ideal (recall that $\B$ is a  valuation ring of a
field $F$ if $\B$ is a local integral domain satisfying $x\in F-\B$ implies $%
x^{-1}\in m$). They called such deformations valued global deformations.

\subsection{Perturbation theory}

The perturbation theory over complex numbers is based on an enlargement of the field of real
numbers with the same algebraic order properties as $\real$. The extension $%
\real^{\star }$ of $\real$ induces the existence of infinitesimal element in $%
\real^{\star }$, hence in $\left( \complex ^{n}\right) ^{\star }$. The '' infinitely
small number'' and ''illimited number'' have a long historical tradition
(Euclid, Eudoxe, Archimedes..., Cavalieri, Galilei...Leibniz, Newton).
The infinitesimal methods was considered in heuristic and intuitive way
until 1960 when A. Robinson gave a rigorous foundation to these methods  
\cite{Robinson60}. He used methods of mathematical logic, he constructed a
NonStandard model for real numbers. However, there exists other frames for
infinitesimal methods(see \cite{Nelson},\cite{Lutz-Goze81}, \cite
{lutz-makhlouf-meyer}). In order to study the local properties of complex Lie algebras 
M. Goze introduced   in 1980 the notion of perturbation of algebraic structure 
(see\cite{Goze88}) in Nelson's framework  \cite
{Nelson}. The description given here is more algebraic, it is based on Robinson's framework.
 First, we summarize a description 
of the hyperreals field and then we define the 
perturbation notion over hypercomplex numbers. 

\subsubsection{Field of hyperreals and their properties}
The construction of hyperreal numbers system needs four
axioms. They induce a triple $(\real,\real^{\star },\star )$
where $\real$ is the real field, $\real^{\star }$ is the hyperreal field and  $\star 
$ is a natural mapping. 

Axiom 1. $\real$ is an archimedean ordered field.

Axiom 2. $\real^{\star }$ is a proper ordered field extension of $\real$.

Axiom 3. For each real function  $f$ of $n$ variables there is a
corresponding  hyperreal function $f^{\star }$ of $n$ variables, called
natural extension of  $f$. The field operations of  $\real^{\star }$ are the
natural extensions of the field operations of $\real$.

Axiom 4. If two systems of formulae have exactly the same real solutions,
they have exactly the same hyperreal solutions. 

The following theorem shows that such extension exists

\begin{theorem}
Let $\real$ be the ordered field of real numbers. There is an ordered field
extension $\real^{\star }$ of $\real$ and a mapping $\star $ from real functions to
hyperreal functions such that the axioms 1-4 hold.
\end{theorem}

In the proof, the plan is to find an infinite set $M$ of formulas $\varphi (x)$ which describe all
properties of a positive infinitesimal $x$, and to built $\real ^{\star}$ out of this set
of formulas. See \cite{Kiesler} pp 23-24 for a complete proof.

\begin{definition}
\begin{itemize}
\item  $x$ $\in \real^{\star }$ is called infinitely small or infinitesimal if $%
\left| x\right| <r$ for all $r\in \real ^{+}$.

\item  $x$ $\in \real^{\star }$ is called limited if there exists $r \in \real$ such that  
$\left| x\right| <r$.

\item  $x$ $\in \real^{\star }$ is called illimited or infinitely large if $\left| x\right| >r$ for
all $r\in \real$

\item  Two elements $x$ and $y$ of $\real^{\star }$ are called infinitely close $%
(x\simeq y)$ if $x-y$ is infinitely small.
\end{itemize}
\end{definition}

\begin{definition}
Given a hyperreal number $x$ $\in \real^{\star }$, we set

$halo(x)=\{y\in \real^{\star },\quad x\simeq y\}$

$galaxy(x)=\{y\in \real^{\star },\quad x-y$ is limited$\}$
\end{definition}
In particular, $halo(0)$ is the set of infinitely small elements of   $%
\real^{\star }$ and $galaxy(0)$ is the set of  limited hyperreal numbers which we denote
by $\real_{L}^{\star }$.

\begin{proposition}
\begin{enumerate}
\item  $\real_{L}^{\star }$ and $halo(0)$ are subrings of $\real^{\star }.$

\item  $halo(0)$ is a maximal ideal in $\real_{L}^{\star }.$
\end{enumerate}
\end{proposition}

\begin{proof}
The first property is easy to prove. (2)  Let $\eta \simeq 0$ and $a\in \real$,
then there exists $t\in \real$ such that $\left| a\right| <t$ and for all $r\in \real$%
, $\left| \eta \right| <r/t$, thus $\left| a\eta \right| <r$. Therefore, $%
a\eta \in \real_{L}^{\star }$ and  $halo(0)$ is an ideal in $\real_{L}^{\star }$ .
Let us show that the ideal is maximal. We note that  $x$ $\in \real^{\star }$ is
illimited is equivalent to $x^{-1}$ is infinitely small. Assume that there
exists an ideal $I$ in $\real_{L}^{\star }$ containing $halo(0)$ and let  $a\in
I\setminus halo(0)$, we have  $a^{-1}\in \real_{L}^{\star }$, thus $1=aa^{-1}\in
I$, so $I\equiv \real_{L}^{\star }.$
\end{proof}
The following theorem shows that there is a homomorphism of  ring $\real_{L}^{\star }$ onto the field of real numbers.

\begin{theorem}
 Every $x$ of $\real_{L}^{\star }$ admits a
unique $x_{0}$ in $\real$ (called standard part of $x$ and denoted by $st(x)
$) such that $x\simeq x_{0}$.
\end{theorem}

\begin{proof}
Let $x\in \real_{L}^{\star }$, assume that it exists two real numbers $r$ and $s$
such that $x\simeq r$ and $x\simeq s$, this implies that $r-s\simeq 0$. Since the unique
infinitesimal real number is $0$ then $r=s$. To show the existence, set
 $X=\left\{ s\in \real:s<x\right\} $. Then $X\neq \emptyset $ and $%
X<r $, where $r$ is a positive real number ($-r<x<r$) . Let  $t$ be the smallest $r$,
for all positive real $r$, we have $x\leq t+r$,
then $x-t\leq r$ and  $t-r<x$, it follows $-(x-t)\leq r$. Then $x-t\simeq 0$, hence $x\simeq t$.
\end{proof}

We check easily that for  $x,y\in \real_{L}^{\star }$,  $st(x+y)=st(x)+st(y)$, 
$st(x-y)=st(x)-st(y)$ and $st(xy)=s(x)st(y)$. Also, for
all  $r\in \real$, $st(r)=r$.

\begin{remark}
All these  concepts and properties may be extended to hyperreal vectors of  $%
\real^{\star ^{n}},$  An element $v=\left( x_{1},\cdots ,x_{n}\right) $, where $%
x_{i}\in \real^{\star }$, is called infinitely small if  $\left| v\right| =\sqrt{%
\sum_{i=1}^{n}x_{i}^{2}}$ is an infinitely small hyperreal and called limited if  $%
\left| v\right| $ is limited. Two vectors $\real^{\star ^{n}}$are infinitely
close if their difference is infinitely small.  The extension to complex
numbers and vectors is similar ($\complex^{\star}=\real^{\star }\times \real^{\star }$).
\end{remark}

\subsubsection{Perturbation of associative algebras.}

We introduce here the perturbation notion of an algebraic structure. Let $%
V=\corps ^{n}$  be 
a $\corps$-vector space ($\corps = \real $ or $\complex$) and
$V^{\star}$ be a vector space over $\corps^{\star}$ 
($\corps^{\star} = \real ^{\star }$ or $\complex^{\star }$). Let $\A$ be an associative
algebra of $Alg_{n}$  with a multiplication $\mu _{0}$  over $V$.

\begin{definition}
A morphism $\mu $ in $Hom((V^{\star})^{\times 2},V^{\star})$ is a
perturbation of $\mu _{0}$  if 
\[
\forall X_{1}, X_{2}\in V\ :\mu \left( X_{1}, 
X_{2}\right) \simeq \mu _{0}\left( X_{1}, X_{2}\right) 
\]
and $\mu$ satisfies the associativity condition over $V$. We write
 $\left( \mu \simeq \mu _{0}\right) $.
\end{definition}

Fixing a basis of $V$, the extension $V^{\star}$ is a vector space which may be taken with the same
basis as $V$. Then $\mu $  is a
perturbation of $\mu _{0}$ is equivalent to say that the difference
between the structure constants of $\mu $ and $\mu _{0}$ is an
infinitesimal vector in $V^{\star }$.

Two perturbations $\mu$ and $\mu^{\prime}$ of $\mu _{0}$ are isomorphic if there is an
invertible map $f$ in $Hom(V^{\star},V^{\star})$ such that $\mu^{\prime}=f^{-1} \circ \mu
\circ (f\otimes f) $.

The following decomposition of a perturbation follows from Goze's
decomposition \cite{Goze88}.

\begin{theorem}
Let $\A$ be an algebra in $ Alg _{n}$ with a multiplication $\mu _{0}$ and let
  $\mu $  be a perturbation of $\mu _{0}$.  We have the following decomposition of $\mu $ :
\[
\mu =\mu _{0}+\varepsilon _{1}\varphi _{1}+\varepsilon _{1}\varepsilon
_{2}\varphi _{2}+....+\varepsilon _{1}\cdot \cdot \cdot \varepsilon
_{k}\varphi _{k}
\]
{\it where}

\begin{enumerate}
\item  $\varepsilon _{1},...,\varepsilon _{k}${\it \ are non zero
infinitesimals in }$\corps ^{\star }.$

\item  $\varphi _{1},...,\varphi _{k}$ are independent bilinear maps
in $Hom(V^{\times 2},V)$.
\end{enumerate}
\end{theorem}

\begin{remark}

\begin{enumerate}
\item  The integer $k$ is called the length of the perturbation. It satisfies 
$k\leq n^{3}$.
\item The perturbation decomposition is generalized in the valued global deformation case 
 \cite{Goze-Remm02} by taking the $\varepsilon _{i}$ in maximal ideal of a
valuation ring.
\item  The associativity of $\mu $ is equivalent to a finite system of
equation called the perturbation equation. This equation is studied above by using
Massey cohomology products .
\end{enumerate}
\end{remark}

\subsubsection{Resolution of the perturbation equation} 
In the following we discuss the conditions on  $\varphi _{1}$ such that it 
is a first term of a perturbation.

Let us consider a perturbation of length 2, $\mu =\mu _{0}+\varepsilon
_{1}\varphi _{1}+\varepsilon _{1}\varepsilon _{2}\varphi _{2},$ the perturbation
equation is equivalent to
\[
\left\{ 
\begin{array}{c}
\delta \varphi _{1}=0 \\ 
\varepsilon _{1}[\varphi _{1}, \varphi _{1}]_{G}+2\varepsilon _{1}\varepsilon
_{2}[\varphi _{1}, \varphi _{2}]_{G}+\varepsilon _{1}\varepsilon
_{2}^{2}[\varphi _{2}, \varphi _{2}]_{G}+2\varepsilon _{2}\delta \varphi _{2}=0
\end{array}
\right. 
\]
Where $[\varphi _{i}, \varphi _{j}]_{G}$ is the trilinear map defined by the 
Gerstenhaber bracket (see section 2.1 ) : 

It follows, as in the deformation equation, that $\varphi _{1}$ is a
cocycle of $Z^{2}\left( \A,\A \right) .$ The second equation has
infinitesimal coefficients but the vectors $$\left\{ [\varphi _{1}, \varphi _{1}]_{G},
[\varphi _{1}, \varphi _{2}]_{G},[\varphi _{2}, \varphi _{2}]_{G},\delta \varphi _{2}\right\} $$ are in 
$Hom(V^{\times
3},V)$ and form a system of rank 1.
\begin{proposition}
{\it Let }$\mu ${\it \ be a perturbation of length }$k${\it \ of the multiplication 
}$\mu _{0}.${\it \ The rank of the vectors }$\left\{ [\varphi _{i},
\varphi _{j}]_{G},\delta \varphi _{i}\right\} \ i=1,..,k$ and $i\leq j\leq k$%
{\it \ is equal to the rank of the vectors }$\left\{ [\varphi _{i},
\varphi _{j}]_{G}\right\} \ i=1,..,k-1$ and $i\leq j\leq k-1.$
\end{proposition}

\begin{proof}
We consider a non trivial linear form $\omega $ containing in its kernel $%
\left\{ [\varphi _{i},
\varphi _{j}]_{G}\right\} \ i=1,..,k-1$ and $i\leq
j\leq k-1$. We apply it to the perturbation equation then it follows that
all the vectors $\left\{ [\varphi _{i},
\varphi _{j}]_{G},\delta \varphi
_{i}\right\} \ i=1,..,k$ and $i\leq j\leq k$ are in the kernel of $\omega $.
\end{proof}
The following theorem, which uses the previous proposition, characterizes the cocycle
 which should be a first term of a perturbation, see \cite{Makhlouf90}
for the proof and \cite{Massey} for the Massey products.
\begin{theorem}
 Let $\A$ be an algebra in $Alg _{n}$ with a multiplication  $\mu _{0}$.
A vector $\varphi _{1}${\it \ in }$Z^{2}\left( \A,\A\right) $ 
{\it is the first term of} {\it a perturbation }$\mu ${\it \ of }$\mu _{0}$%
( $k$ {\it being the length of }$\mu $) if and only if

\begin{enumerate}
\item  {\it The massey products }$\left[ \varphi _{1}^{2}\right] ,\left[
\varphi _{1}^{3}\right] \ldots \left[ \varphi _{1}^{p}\right] ${\it \ vanish
until }$p=k^{2}.$

\item  {\it The product representatives in }$B^{3}\left( \A,\A \right) ${\it \ form a system of rank less or equal to }$k\left(
k-1\right) /2.$
\end{enumerate}
\end{theorem}

\subsection{Coalgebra perturbations}
The perturbation notion can be generalized to any algebraic structure on $V$. 
A morphism $\mu $ in 
$Hom((V^{\star})^{\times p},(V^{\star})^{\times q})${\it \ is a
perturbation of }$\mu _{0}$ in 
$Hom(V^{\times p},V^{\times q})${\it \ }$\left( \mu \simeq \mu _{0}\right) ${\it if} 
\[
\forall X_{1},\cdots X_{p}\in V\ :\mu \left( X_{1}, \cdots ,
X_{p}\right) \simeq \mu _{0}\left( X_{1}, \cdots , X_{p}\right) 
\]

In particular, if $p=1$ and $q=2$, we have the concept of coalgebra perturbation. 
\section{Comparison of the deformations}
We show that the formal deformation and the perturbation are global deformations with 
appropriate bases, and we show that over complex numbers the perturbations contains all
the convergent deformations.
\paragraph{Global deformation and formal deformation}
The following proposition gives the link between formal deformation and global deformation.
\begin{proposition}
Every formal deformation is a global deformation.
\end{proposition}
Every formal deformation of an algebra $\mathcal{A}$, in Gerstenhaber sense, is a global deformation with
a basis $\left( \B,m\right) $ where $\B=\corps[[t]]$ and $m=t \corps[[t]]$.

\paragraph{Global deformation and perturbation}
Let $\corps  ^{\star }$ be a proper extension of $\corps$ described in the section (2.4)
, 
where $\corps$ is the field  $\real$ or $\complex$. 
\begin{proposition}
Every perturbation of an algebra over $\corps$  is a global 
deformation with a 
base $ \left( \corps _{L}^{\star },halo(0)\right) $. 
\end{proposition}
\begin{proof}
The set $\B:= \corps _{L}^{\star }$ is a local ring
formed by the limited elements of $ \corps  ^{\star }$. We
define the augmentation $\varepsilon : \corps _{L}^{\star
}\rightarrow \corps$ which associates to any $x$ its standard part $st(x)$.
The Kernel of $\varepsilon $ corresponds to $halo(0)$,  the set of infinitesimal
elements of $ \corps ^{\star }$, which is a maximal ideal in $\corps _{L}^{\star
}$.
\end{proof}
The multiplication $\mu _{\B }$ of a global deformation with a base $%
\left(  \corps _{L}^{\star },halo(0)\right) $ of an algebra with
a multiplication $\mu $ may be written, for all $x,y\in \corps^{n}$ : 
\[
\mu _{B }(1\otimes x,1\otimes y)=1\otimes \mu (x,y)+\sum_{i}\alpha
_{i}\otimes z_{i}\quad \text{avec }\alpha _{i}\in halo(0)\subset \corps^{\star},\ z_{i}\in \corps^{n} 
\]
This is equivalent to say that $\mu _{\B}$ is a perturbation of $\mu $.

\begin{remark}
The corresponding global deformations of a perturbations are local but not
infinitesimal because $halo(0)^{n}\neq 0$.
\end{remark}
\paragraph{Formal deformation and perturbation}
\begin{proposition}
Let $\A$ be  a $\complex$-algebra in $Alg_n$ with a multiplication $\mu_{0}$. Let $\mu_{t}$
be a convergent deformation of $\mu_{0}$ and $\alpha$ be infinitesimal in $\complex^{\star}$ 
then $\mu_{\alpha}$ is a perturbation of $\mu_{0}$.
\end{proposition}
\begin{proof}
Since the deformation is convergent then  $\mu_{\alpha}$ corresponds to a point of
$Alg_n \subset \complex ^{n^3}$ which is infinitely close to the corresponding point of $\mu_{0}$. Then it determines a
perturbation of $\mu_{0}$.

\end{proof}

A perturbation should correspond to a formal deformation if we consider the power series ring
$\corps[[t_{1},\cdots, t_{r}]]$ with more than one parameter. 

\section{Universal and versal deformation}

Given an algebra, the problem is to find particular deformations which induces all the others in 
the space of all deformations or in a fixed
category of deformations. We say that the deformation is universal if there is unicity of the
homomorphism between base algebras. This problem is too hard in general but the global deformation seems 
more adapted to construct universal or versal deformation.   This problem was considered in Lie
algebras case for the categories of deformations over infinitesimal local algebras and complete 
local algebras  (see
 \cite{Fialowski88},\cite{Fialowski99},
 \cite{Fialowski02}). They show that if we consider the
infinitesimal deformations, i.e. the deformations over local algebras $B$ such that
$m_{B}^{2}=0$ where  $m_{B}$ is tha maximal ideal, then there exists a universal
deformation. If we consider the category of complete local rings, then there does
not exist a universal deformation but only versal deformation.
 
\paragraph{Formal global deformation.}

Let $\B$ be a complete local algebra over $\corps$, so $\B=$ 
$\overleftarrow{\lim }_{n\rightarrow \infty }(\B/m^{n})$  (inductive limit), where $m$ is the
maximal ideal of $\B$ and we assume that $\B/m\cong \corps$.

\begin{definition}
A formal global deformation of $\mathcal{A}$ with base $\left( \B,m\right) $ is an
algebra structure on the completed tensor product $\B\stackrel{\wedge }{%
\otimes }\mathcal{A} =\overleftarrow{\lim }_{n\rightarrow \infty }
((\B/m^{n})\otimes \mathcal{A} )$ such that $\varepsilon \stackrel{\wedge }{%
\otimes }id:\B\stackrel{\wedge }{\otimes }\mathcal{A} \rightarrow K\otimes 
\mathcal{A} =\mathcal{A} $ is an algebra homomorphism.
\end{definition}

\begin{remark}
\begin{itemize}
\item The formal global deformation of $\mathcal{A} $ with base $\left(
\corps [[t]],t\corps [[t]]\right) $ are the same as formal one parameter deformation
of Gerstenhaber.
\item The perturbation are complete global deformation because
$lim_{n\rightarrow \infty }\corps^{\star}/halo(0)^{n}$ is isomorphic to $\corps^{\star}$ and $\corps^{\star}$ is
isomorphic $\corps$ as algebra.
\end{itemize}
\end{remark}

We assume now that the algebra $\A$ satisfies $\dim \left( H^{2}\left( \mathcal{A}
,\mathcal{A} \right) \right) <\infty $. We consider $\B=\corps \oplus H^{2}\left(
\mathcal{A} ,\mathcal{A} \right) ^{dual }$. The following theorems du to
Fialowski and Post \cite{Fialowski88}\cite{Fialowski99} show the existence of universal infinitesimal deformation under the previous
assumptions.

\begin{theorem}
There exists, in the category of infinitesimal global deformations, a universal infinitesimal deformation $\eta _{%
\mathcal{A} }$ with base $\B$ equipped with the multiplication $\left( \alpha
_{1},h_{1}\right) \cdot \left( \alpha _{1},h_{1}\right) =\left( \alpha
_{1}\alpha _{2},\alpha _{1}h_{2}+\alpha _{2}h_{1}\right) $. 
\end{theorem}
Let $\mathbb{P}$ be any finitedimensional local algebra over $\corps$. The theorem means that
for any infinitesimal deformation defined by $\mu_{\mathbb{P}} $ of an algebra $\mathcal{A} $,   
there exists a unique homomorphism $\Phi
:\corps\oplus H^{2}\left( \mathcal{A} ,\mathcal{A} \right) ^{dual }\rightarrow \mathbb{P}$ such
that $\mu_{\mathbb{P}} $ is equivalent to the push-out $\Phi _{*}\eta _{\mathcal{A} }$.
\begin{definition}
A formal global deformation $\eta $ of $\mathcal{A} $ parameterized by a complete local
algebra $\B$ is called versal if for any deformation $\lambda $ of $\mathcal{A}$,
parameterized by a complete local algebra $\left( \mathbb{A},m_{\mathbb{A}}\right) $, there is
a morphism $f:\B\rightarrow \mathbb{A}$ such that \newline
(1) The push-out $f_{*}\eta $ is equivalent to $\lambda $.\newline
(2) If $\mathbb{A}$ satisfies $m_{\mathbb{A}}^{2}=0$ , then $f$ is unique;
\end{definition}

\begin{theorem}
Let $\A$ be an algebra .

(1) There exists a versal formal global deformation of $\mathcal{A} $.
\newline
(2) The base of the versal formal deformation is formally embedded into $
H^{2}\left( \mathcal{A} ,\mathcal{A} \right) $ (it can be described in $H^{2}\left( 
\mathcal{A} ,\mathcal{A} \right) $ by a finite system of formal equations)
\end{theorem}

\paragraph{Examples.}
The Witt algebra is the infinite dimensional Lie algebra of polynomial vector fields spanned by
the fields $e_{i}=z^{i+1}\frac{d}{dz}$ with $i \in \mathbb{Z}$.
In \cite{Fialowski88}, Fialowski constructed versal deformation of
the Lie subalgebra $L_{1}$ of Witt algebra ($L_{1}$ is spanned by $e_{n}$, $n>0$, while the bracket
is given by $[e_{n},e_{m}]=(m-n)e_{n+m}$ ). 

Here is three real deformations of the Lie algebra $L_1$ which are non trivial
and pairwise nonisomorphic.

$[ e_i, e_j ]_{t}^1 = ( j - i ) ( e_{i + j} + te_{i + j - 1}$

$[ e_i, e_j ]_{t}^2 = \{ 
\begin{array}{l}
  ( j - i ) e_{i + j}  \text{ if } i, j > 1\\
  ( j - i ) e_{i + j} + tj e_j  \text{ if} i = 1
\end{array}$

$[ e_i, e_j ]_{t}^3 = \{ 
\begin{array}{l}
  ( j - i ) e_{i + j}  \text{ if} i, j \neq 2\\
  ( j - i ) e_{i + j} + tje_j  \text{ if } i = 2
\end{array}$

Fialowski and Post, in \cite{Fialowski99},  study the
$L_{2}$ case. A more general procedure using the  Harrisson cohomology of the commutative algebra
$\B$ is described by Fialowski and Fuchs
in \cite{Fialowski00}.

The Fialowski's global deformation of $L_{1}$ can be realized as a
perturbation, the parameter $t_{1},t_{2\text{ }}$ and $t_{3}$ have to be different 
infinitesimals in $\complex ^{\star }$. 

\section{Degenerations}
The notion of degeneration is fundamental in the geometric study of $Alg_{n}$ and helps,
in general, to construct new algebras.
\begin{definition}
Let $\A _{0}$ and $\A_{1}$ be two n-dimensional algebras. We say  that $\A _{0}$ is a  degeneration
 of $\A _{1}$ if $\A_{0}$ belongs to $\overline{\vartheta \left( \A_{1}\right) }$ , 
 the Zariski closure of the orbit of $\A_{1}$.
 \end{definition}
In the following we define degenerations in the different frameworks.
 \subsection{Global viewpoint}
A characterization of global degeneration  was given by Grunewald and O'Halloran in
 \cite{Grun_OHallo} :
 \begin{theorem}
 Let $\mathcal{A}_0$ and $\mathcal{A}_1$ be two $n$-dimensional associative
algebras over $\corps$ with the multiplications $\mu_0$ and $\mu_1$. The
algebra $\mathcal{A}_0$ is a degeneration of $\mathcal{A}_1$ if and only if
there is a discrete valuation $\mathbb{K}$-algebra $\B$ with residue field $\corps$ whose quotient field
$\mathcal{K}$ is
finitely generated over $\corps$ of transcendence degree one (one parameter), and there is an $n$-dimensional algebra $\mu
_{\B}$ over $\B$ such that
 $\mu_{\B} \otimes \mathcal{K} \cong \mu_1 \otimes \mathcal{K}$ and $\mu_{\B} \otimes \corps \cong \mu_0$.

 \end{theorem}
 \subsection{Formal viewpoint}
Let $t$ be a parameter in $\corps$ and $\{f_{t}\}_{t\neq 0}$\ be a continuous 
family of  invertible linear maps on $V$\ over $\corps$\ and $\A_{1} =(V,\mu_{1}) $ be an 
algebra over $\corps$. 
The limit (when it exists) of a sequence $f_{t}\cdot \A_{1}$\ , $%
\A _{0}=\lim_{t\rightarrow 0}f_{t}\cdot \A_{1}$ , is a {\em formal degeneration \em} of $\A_{1}$ in the
sense that $\A _{0}$\ is in the Zariski closure of the set $\left\{ f_{t}\cdot
\A\right\} _{t\neq 0}.$
\\  The multiplication $\mu _{0}$ is given by 
$$
\mu _{0}=\lim_{t\rightarrow 0}f_{t}\cdot \mu_{1} =\lim_{t\rightarrow
0}f_{t}^{-1}\circ \mu_{1} \circ f_{t}\times f_{t} 
$$
\begin{itemize}
\item The multiplication $\mu _{t}=f_{t}^{-1}\circ \mu_{1} \circ f_{t}\times f_{t}$
satisfies the associativity condition. Thus, when $t$
tends to $0$ the condition remains satisfied.

\item  The linear map $f_{t}$\ is invertible when $t\neq 0$ \ and may be
singular when $t=0.$ Then, we may obtain by degeneration a new algebra.

\item The definition of formal degeneration maybe extended naturally to infinite dimensional case.

\item  When $\corps$ is the complex field, the multiplication given by the limit, follows 
from a limit of the structure constants,
using the metric topology. In fact, $f_{t}\cdot \mu $ 
corresponds to a change of basis when $t\neq 0$. When $t=0$, they give eventually a new
point in $Alg_{n}\subset \corps^{n^{3}}$.

\item  If $f_t$ is defined by a power serie the images over $V\times V$ of the multiplication 
of $%
f_{t}\cdot \A$ are in general in the Laurent power series ring $V\left[
\left[ t,t^{-1}\right] \right] $.  But when the
degeneration exists, it lies in the power series ring $V\left[ \left[ t\right]
\right] $. 
\end{itemize}
\begin{proposition}
Every formal degeneration is a global degeneration.
\end{proposition}
The proof follows from the theorem (5.1) and the last remark.

\subsection{Contraction}
The  notion of degeneration over the hypercomplex field is called contraction. It is defined by :
\begin{definition}
Let $\A_{0}$ and $\A_{1}$ be two
algebras in $Alg _{n}$ with multiplications $\mu _{0}$ and $\mu _{1}$. 
The algebra $\A _{0}$ is a contraction of $\A_{1}$ if there exists a perturbation 
$\mu $ of $\mu _{0}$ such that $\mu $ is
isomorphic to $\mu _{1}$.
\end{definition}
 The definition gives a characterization over hypercomplex field of $\mu _{0}$  in the closure of the orbit of $\mu _{1}$ (for the
usual topology of $\complex^{n^{3}}$).

\subsection{Examples}

\begin{enumerate}
\item {\it The null algebra of }$Alg _{n}$ {\it is a degeneration of every algebra
of }$Alg _{n}$.\newline
In fact, the null algebra is given in a basis $\left\{
e_{1},....,e_{n}\right\} $ by the following nontrivial products $\mu _{0}\left( e_{1},e_{i}\right) =\mu_{0} \left(
e_{i},e_{1}\right) =e_{i}\ \quad i=1,...,n.$ \newline
Let $\mu_{1} $ be a multiplication of any algebra of $Alg _{n}$, we have $\mu
_{0}=\lim_{t\rightarrow 0}f_{t}\cdot \mu_{1} $ with $f_{t}$ given by the
diagonal matrix $\left( 1,t,....,t\right)$. \newline
For $i\neq 1$ and $ j\neq 1$ we have 
$$
f_{t}\cdot \mu \left( e_{i},e_{j}\right) =f_{t}^{-1}\left( \mu \left(
f_{t}\left( e_{i}\right) ,f_{t}\left( e_{j}\right) \right) \right)
=f_{t}^{-1}\left( \mu \left( te_{i},te_{j}\right) \right) = 
$$
$$t^{2}C_{ij}^{1} e_{1}+t \sum_{k > 1} C_{ij}^k e_k \rightarrow 0 (\text{when }t 
\rightarrow 0)
$$
For $i=1$ and $j\neq 1$ we have 
\[
f_{t}\cdot \mu \left( e_{1},e_{j}\right) =f_{t}^{-1}\left( \mu \left(
f_{t}\left( e_{1}\right) ,f_{t}\left( e_{j}\right) \right) \right)
=f_{t}^{-1}\left( \mu \left( e_{1},te_{j}\right) \right) =f_{t}^{-1}\left(
te_{j}\right) =e_{j}.
\]
This shows that every algebra of $Alg_{n}$ degenerates formally to a null algebra.

By taking the parameter $t=\alpha$, $\alpha$ infinitesimal in the field of hypercomplex numbers, we get that the null algebra
is a contraction of any complex associative algebra in $Alg_{n}$. In fact, $f_{\alpha}\cdot \mu_{1}$ is a
perturbation of $\mu_{0}$ and is isomorphic to $\mu_{1}$.

The same holds also in the global viewpoint : Let $B=\corps [t]_{<t>}$ be the polynomial ring localized at the
prime ideal $<t>$ and let $\mu_{B}=t \mu_{1}$ on the elements different from the unity and $\mu_{B}=\mu_{1}$
elsewhere. Then $\mu_{1}$ is $\corps (t)$-isomorphic to $\mu_{B}$ via $f_{t}^{-1}$ and $\mu_{B}\otimes
\corps=\mu_{0}$
\item  The classification of $Alg _{2}$ yields two isomorphic classes.

Let $\left\{ e_1,e_2\right\} $ be a basis of $\corps ^2.$

$\A_1\ :\ \mu _1\left( e_1,e_i\right) =\mu _1\left( e_i,e_1\right) =e_i\
\quad i=1,2;\ \mu _1\left( e_2,e_2\right) =e_2.$

$\A_0\ :\ \mu _0\left( e_1,e_i\right) =\mu _0\left( e_i,e_1\right) =e_i\
\quad i=1,2;\ \mu _0\left( e_2,e_2\right) =0.$

Consider the formal deformation $\mu _t$ of $\mu _0$ defined by

\[
\mu _t\left( e_1,e_i\right) =\mu _t\left( e_i,e_1\right) =e_i\ \quad i=1,2;\
\mu _t\left( e_2,e_2\right) =t\ e_2. 
\]

Then $\mu _t$ is isomorphic to $\mu _1$ through the change of basis given by
the matrix $f_t=\left( 
\begin{array}{cc}
1 & 0 \\ 
0 & t
\end{array}
\right) $,  $\mu _t=f_t\cdot \mu _1.$

We have $\mu _{0}=$ $\lim_{t\rightarrow 0}f_{t}\cdot \mu _{1},$ thus $\mu
_{0}$ is a degeneration of $\mu _{1}$. Since $\mu _{1}$ is isomorphic to $\mu
_{t}$, it is a deformation of $\mu _{0}$.
We can obtain the same result over complex numbers if we consider the parameter $t$
infinitesimal in $\complex ^{\star}$.

\end{enumerate}

\subsubsection{Graded algebra}
In this Section, we give a relation between an algebra and its associated graded
algebra.
\begin{theorem}
Let \( \mathcal{A} \) be an algebra over $\corps$ and \( \A _{0} \)
\( \subseteq \A _{1}\subseteq \cdots \subseteq \A _{n}\subseteq \cdots  \)
an algebra filtration of \( \A \), \( \A =\cup _{n\geq 0}\A _{n} \). Then the
graded algebra \( gr(\A)=\oplus _{n\geq 1}\A_{n}/\A_{n-1} \) is
a formal degeneration of \( \A \).
\end{theorem}
\begin{proof}
Let \( \A [[t]] \) be a power series ring in one variable $t\in \corps$ over \( \A \), 
\( \A[[t]]=\A\otimes \corps[[t]]=\oplus _{n\geq 0}\A \otimes t^{n} \).

We denote by \( \A_{t} \) the Rees algebra associated to the filtered
algebra \( \A \), \( \A_{t}=\sum _{n\geq 0}\A_{n}\otimes t^{n} \).
The Rees algebra \( \A_{t} \) is contained in the algebra \( \A[[t]] \).

For every \( \lambda \in \corps \), we set \( \A_{(\lambda )}=\A_{t}/((t-\lambda )\cdot \A_{t}) \).
For \( \lambda =0 \), \( \A_{(0)}=\A_{t}/(t\cdot \A_{t}) \). The algebra 
\( \A_{(0)} \) corresponds to the graded algebra $gr(\A)$ and $\A_{(1)}$ is
isomorphic to $\A$. In fact, we suppose that the parameter
\( t \) commutes with the elements of \( \A \) then
 \( t\cdot \A_{t}= \)\( \oplus _{n\geq 0}\A_{n}\otimes t^{n+1} \). 
It follows that $
\A_{(0)}=\A_{t}/(t\cdot \A_{t})=
(\sum _{n}\A_{n}\otimes t^{n})/(\oplus _{n}\A_{n}\otimes t^{n+1})=
\oplus _{n\geq 0}\A_{n}/\A_{n-1}=gr(\A)$.
By using the linear map from \( \A_{t} \) to \( \A \) where the image
of \( a_{n}\otimes t^{n} \) is \( a_{n} \), we have \( \A_{(1)}=\A_{t}/((t-1)\cdot \A_{t})\cong \A \).

If \( \lambda \neq 0 \),
the change of parameter \( t=\lambda T \) shows that \( \A_{(\lambda )} \)
is isomorphic to \( \A_{(1)} \). This ends the proof that \( \A_{(0)}=gr(\A) \)
is a degeneration of \( \A_{(1)}\cong \A \). 
\end{proof}
\subsection{Connection between degeneration and deformation}
The following proposition gives a connection between degeneration and
deformation. 

\begin{proposition}
Let $\A_{0}$ and $\A_{1}$ be two algebras in $Alg_n$.
If $\A _{0}$ is a degeneration of $\A _{1}$ then $\A _{1}$
 is a deformation of $\A _{0}.$
\end{proposition}

In fact, let $\A _{0}=$ $\lim_{t\rightarrow 0}f_{t}\cdot \A_{1}$ be a formal degeneration
of $\A$ then $\A _{t}=f_{t}\cdot \A_{1}$ is a formal deformation of $\A _{0}$.

In the contraction sense such a property follows directly from the definition.

In the global viewpoint, we get also that every degeneration can be realized by a global deformation. The base of
the deformation is the completion of the discrete valuation $\corps$-algebra (inductive limit of
$\mu_{n}=\mu_{B} \otimes B/m_{B}^{n+1}$). 

\begin{remark}
The converse is, in general, false. The following example shows that there is no
duality in general between deformation and degeneration.
\end{remark}
We consider, in $Alg _{4}$, the family $\A _{t}=\complex \left\{ x,y\right\}
/\left\langle x^{2},y^{2},yx-txy\right\rangle $ where $\complex \left\{ x,y\right\} $
stands for the free associative algebra with unity and generated by $x$ and $y$.

Two algebras $\A _{t}$ and $\A _{s}$ with $t\cdot s\neq 1$ are not isomorphic. They are
isomorphic when $t \cdot s=1$.
\newline
Thus, $\A _{t}$ is a deformation of $\A _{0}$ but the family $\A _{t}$ is not
isomorphic to one algebra and cannot be written $\A _{t}=f_{t} \cdot \A _{1}$.

We have also more geometrically the following proposition :
\begin{proposition}
If an algebra $\A_{0}$ is in the boundary of the orbit of $\A_{1}$, then this degeneration defines a
nontrivial deformation of $\A_{0}$.
\end{proposition}

 With the contraction viewpoint, we can  characterize the perturbation arising from degeneration by :
 \begin{proposition}
 Let  $\mu $ be a perturbation over hypercomplex numbers of an algebra with a multiplication $\mu
_{1}$. Then $\mu $ arises from a contraction if there exists an
element $\mu _{0}$ with structure constants in $\mathbb{C}$ belonging to the orbit of $\mu$.
 
\end{proposition}
In fact, if $\mu$ is a perturbation of $\mu _0$ and there exists an algebra with
multipication $\mu_1$ such that: $\mu \in \vartheta \left( \mu _1\right) $ and $\mu
\simeq \mu _0.$ Thus $\mu _0$ is a contraction of $\mu _1.$ This shows that the orbit of $\mu$ passes
through a point of $Alg_{n}$ over $\complex$.

\subsection{Degenerations with $f_{t}=v+t\cdot w$}
Let $f_{t}=v+tw$ be a family of endomorphisms where $v$ is a singular linear
map and $w$ is a regular linear map. The aim of this section is to find a
necessary and sufficient conditions on $v$ and $w$ such that a degeneration
of a given algebra $\A=\left( V,\mu \right) $ exists$.$ We can set $w=id$
because $f_{t}=v+tw=\left( v\circ w^{-1}+t\right) \circ w$ which is
isomorphic to $v\circ w^{-1}+t$. Then with no loss of generality we can
consider the family $f_{t}=\varphi +t\cdot id$ from $V$ into $V$ where $%
\varphi $ is a singular map. The vector space $V$ can be decomposed by $%
\varphi $ under the form $V_{R}\oplus V_{N}$ where $V_{R}$ and $V_{N}$ are $%
\varphi $-invariant defined in a canonical way such that $\varphi $ is
surjective on $V_{R}$ and nilpotent on $V_{N}$. Let $q$ be the smallest
integer such that $\varphi ^{q}\left( V_{N}\right) =0$. The inverse of $%
f_{t} $ exists on $V_{R}$ and is equal to $\varphi ^{-1}\left( t\varphi
^{-1}+id\right) ^{-1}$. But on $V_{N},$ since $\varphi ^{q}=0$, it is given
by 
\[
\frac{1}{\varphi +t\cdot id}=\frac{1}{t}\cdot \frac{1}{\varphi /t+id}=\frac{1%
}{t}\cdot \sum_{i=0}^{\infty }\left( -\frac{\varphi }{t}\right) ^{i}=\frac{1%
}{t}\cdot \sum_{i=0}^{q-1}\left( -\frac{\varphi }{t}\right) ^{i} 
\]

It follows that 
\[
f_t^{-1}=\left\{ 
\begin{array}{c}
\varphi ^{-1}\left( t\varphi ^{-1}+id\right) ^{-1}\quad \text{on }V_R \\ 
\frac 1t\cdot \sum_{i=0}^{q-1}\left( -\frac \varphi t\right) ^i\quad \text{%
on }V_N
\end{array}
\right. 
\]

Let $f_t=\varphi +t\cdot id$ be a family of linear maps on $V$, where $%
\varphi $ is a singular map. The action of $f_t$ on $\mu $ is defined by $%
f_t\cdot \mu =f_t^{-1}\circ \mu \circ f_t \times f_t$ then 
\[
\begin{array}{c}
f_t\cdot \mu \left( x,y\right) =f_t^{-1}\circ \mu \left( f_t\left( x\right)
,f_t\left( y\right) \right) \\ 
=f_t^{-1}\left( \mu \left( \varphi \left( x\right) ,\varphi \left( y\right)
\right) +t\left( \mu \left( \varphi \left( x\right) ,y\right) +\mu \left(
x,\varphi \left( y\right) \right) \right) +t^2\mu \left( x,y\right) \right)
\end{array}
\]
Since every element $v$ of $V$ decomposes in $v=v_R+v_N$, we set

\[
\begin{array}{c}
A=\mu \left( x,y\right) =A_R+A_N, \\ 
B=\mu \left( \varphi \left( x\right) ,y\right) +\mu \left( x,\varphi \left(
y\right) \right) =B_R+B_N \\ 
C=\mu \left( \varphi \left( x\right) ,\varphi \left( y\right) \right)
=C_R+C_N
\end{array}
\]

Then 
\[
f_t\cdot \mu \left( x,y\right) =\varphi ^{-1}\left( t\varphi ^{-1}+id\right)
^{-1}\left( t^2A_R+tB_R+C_R\right) +\frac 1t\cdot \sum_{i=0}^{q-1}\left(
-\frac \varphi t\right) ^i\left( t^2A_N+tB_N+C_N\right) 
\]

If the parameter $t$ goes to $0,$ then $\varphi ^{-1}\left( t\varphi
^{-1}+id\right) ^{-1}\left( t^2A_R+tB_R+C_R\right) $ goes to $\varphi
^{-1}\left( C_R\right) .$ The limit of the second term is :

\begin{center}
$\lim _{t\rightarrow 0}\frac 1t\cdot \sum_{i=0}^{q-1}\left( -\frac \varphi
t\right) ^i\left( t^2A_N+tB_N+C_N\right) $

$=\lim _{t\rightarrow 0}\left( \frac{id}t-\frac \varphi {t^2}+\frac{\varphi
^2}{t^3}-\cdots \left( -1\right) ^{q-1}\frac{\varphi ^{q-1}}{t^q}\right)
\left( t^2A_N+tB_N+C_N\right) $

$=\lim _{t\rightarrow 0}tA_N-\varphi \left( A_N\right) +\left( \frac{id}%
t-\frac \varphi {t^2}+\frac{\varphi ^2}{t^3}-\cdots \left( -1\right) ^{q-1}%
\frac{\varphi ^{q-3}}{t^{q-2}}\right) \varphi ^2\left( A_N\right) +$

$B_N-\left( \frac{id}t-\frac \varphi {t^2}+\frac{\varphi ^2}{t^3}-\cdots
\left( -1\right) ^{q-1}\frac{\varphi ^{q-2}}{t^{q-1}}\right) \varphi \left(
B_N\right) +$

$\left( \frac{id}t-\frac \varphi {t^2}+\frac{\varphi ^2}{t^3}-\cdots \left(
-1\right) ^{q-1}\frac{\varphi ^{q-1}}{t^q}\right) C_N$

$=\lim _{t\rightarrow 0}tA_N+B_N-\varphi \left( A_N\right) +\qquad \qquad
\qquad $

$\left( \frac{id}t-\frac \varphi {t^2}+\frac{\varphi ^2}{t^3}-\cdots \left(
-1\right) ^{q-1}\frac{\varphi ^{q-1}}{t^q}\right) \left( \varphi ^2\left(
A_N\right) -\varphi \left( B_N\right) +C_N\right) $
\end{center}

This limit exists if and only if 
\[
\varphi ^2\left( A_N\right) -\varphi \left( B_N\right) +C_N=0 
\]
which is equivalent to 
\[
\quad \varphi ^2\left( \mu \left( x,y\right) _N\right)
-\varphi \left( \mu \left( \varphi \left( x\right) ,y\right) _N\right)
-\varphi \left( \mu \left( x,\varphi \left( y\right) \right) _N\right) +\mu
\left( \varphi \left( x\right) ,\varphi \left( y\right) \right) _N=0 
\]
And the limit is $B_N-\varphi \left( A_N\right) .$

\begin{proposition}
 The degeneration of an algebra with a multiplication $\mu $ exists if and
only if the condition 
\[
\varphi ^{2}\circ \mu _{N}-\varphi \circ \mu
_{N}\circ (\varphi \times id)-\varphi \circ \mu _{N}\circ (id\times
\varphi )+\mu _{N}\circ (\varphi \times \varphi )=0
\]
\ where $\mu _{N}\left( x,y\right) =\left( \mu \left( x,y\right) \right)
_{N}$,  holds. And it is defined by 
\[
\mu _{0}=\varphi ^{-1}\circ \mu _{R}\circ (\varphi \times \varphi )+\mu
_{N}\circ (\varphi \times id)+\mu _{N}\circ (id\times \varphi )-\varphi
\circ \mu _{N}
\]
\end{proposition}
 
\begin{remark}
Using a more algebraic definition of a degeneration, Fialowski and O'Halloran \cite{Fialowski90} show that we have
also a notion of universal degeneration and such degeneration exists for every finite dimensional
Lie algebra.  The definition and the result holds naturally in the associative algebra case.
\end{remark}

\section{Rigidity}
An algebra which has no nontrivial deformations is called rigid. In the finite dimensional case, 
this notion is related geometrically
to open orbits in $Alg_{n}$. The Zariski closure of open orbit determines an irreducible component of
$Alg_{n}$.  

\begin{definition}
An algebra $\A$  is called {\em algebraically rigid \em} if $H^2\left(\A,\A\right) =0$.

An algebra $\A$ is called {\em formally rigid \em} if every formal deformation
of $\A$ is trivial.

An algebra $\A$ is called {\em geometrically rigid \em} if its orbit is Zariski open.
\end{definition}
We have the following equivalence du to Gerstenhaber and Schack \cite{Gerstenhaber-Schack}.
\begin{proposition}
Let $\A$ be a finite-dimensional algebra
over a field of characteristic $0$. Then formal  rigidity of $\A$ is equivalent to geometric
rigidity.
\end{proposition}

As seen in theorem (2.1), $H^2\left(
\A,\A\right) =0$ implies that every formal deformation is trivial, then algebraic rigidity 
implies formal rigidity. But the converse is
false. An example of algebra which is formally rigid but not algebraically rigid, {\it $H^2\left(
\A,\A\right) \neq 0$ }, was given by Gerstenhaber and Schack in positive
characteristic and high dimension \cite{Gerstenhaber-Schack}. We do not know such examples in
characteristic $0$. However, there is many rigid Lie algebras, in characteristic $0$, with a nontrivial
 second Chevalley-Eilenberg cohomology group.
 
 Now, we give a sufficient condition for the formal rigidity of an algebra $\A$ using 
the following map $$Sq : \begin{array}{l}
  H^2 ( \A, \A ) \rightarrow H^3 ( \A, \A )\\
  \mu_1 \rightarrow Sq ( \mu_1 ) = \mu_1 \circ \mu_1
\end{array}$$

Let $\mu_1 \in Z^2 ( \A, \A )$, $\mu_1$ is integrable if $Sq ( \mu_1 ) =
0$. If we suppose that $Sq$ is injective then $Sq ( \mu_1 ) = 0$
implies that the cohomology class $\mu_1 = 0$. Then every integrable
infinitesimal is equivalent to the trivial cohomology class. Therefore every formal 
deformation is trivial.

\begin{proposition}
If the map $Sq $ is injective then A is formally rigid.
\end{proposition}

We have the following definition for rigid complex algebra with the perturbation viewpoint.
\begin{definition}
A complex algebra of $Alg_{n}$ with multiplication $\mu _0$ is {\em infinitesimally rigid \em} if every
perturbation $\mu $ of $\mu _0$ belongs to the orbit $\vartheta \left(
\mu _0\right)$. 
\end{definition}
This definition characterizes the open sets over hypercomplex field (with metric topology). 
Since open orbit (with metric topology) is Zariski open \cite{Carles78}. Then
\begin{proposition}
For finite dimensional complex algebras, infinitesimal  rigidity is equivalent to geometric rigidity,
thus to formal rigidity.
\end{proposition}
In the global viewpoint we set the following two concepts of rigidity.
\begin{definition}
Let $\B$ be a commutative algebra over a field $\corps$ and $m$ be a maximal ideal of $\B$.
An algebra $\A$ is called {\em $(\B,m)$-rigid \em}if every
global deformation parametrized by $(\B,m)$ is isomorphic to $\A$ (in the push-out sense).

An algebra is called {\em globally rigid \em} if for every commutative algebra $\B$ and a
maximal ideal $m$ of $\B$ it is $(\B,m)$-rigid
\end{definition}

The global rigidity implies the formal rigidity but the converse is false. Fialowski and
Schlichenmaier show that over complex field the Witt algebra which is algebraically and formally
rigid is not globally rigid. They use families of Krichever-Novikov type algebras
\cite{Fialowski02}.

Finally, we summarize the link between the diffrent concepts of rigidity in the following theorem.

\begin{theorem}
Let $\A$ be a finite dimensional algebra over a field $\corps$ of characteristic zero. Then
\begin{center}
algebraic rigidity $\Longrightarrow$ formal rigidity

formal rigidity $\Longleftrightarrow$ geometric rigidity $\Longleftrightarrow$
$(\corps[[t]],t\corps[[t]])$-rigidity 

global rigidity $\Longrightarrow$ formal rigidity
\end{center}
In particular, if $\corps = \complex$
\begin{center}
infinitesimal rigidity $\Longleftrightarrow$ formal rigidity
\end{center}
\end{theorem}

Recall that semisimple algebras are algebraically rigid. They are classified by Wedderburn's theorem. The classification of low
dimensional rigid algebras is known until $n<7$, see \cite{Gabriel74},\cite{Mazzola79} and 
\cite{Makhlouf-Goze96}. 
The classification of $6$-dimensional rigid algebras  was given using perturbation methods in
\cite{Makhlouf-Goze96}.

\section{The algebraic varieties $Alg_{n}$ }
A point in $Alg_{n}$ is defined by $n^3$ parameters, which are the structure constants
$C_{ij}^{k}$, satisfying a finite system of quadratic relations given by the
associativity condition.  
The orbits  are in 1-1-correspondence with the isomorphism
classes of $n$-dimensional algebras. 
The stabilizer subgroup of $\A$ $\left( stab\left( \A\right) =\left\{ f\in
GL_n\left( \corps\right) :\A=f\cdot \A\right\} \right) $ is  $Aut\left(
\A\right) $, the automorphism group of $\A$. 
The orbit $\vartheta \left( \A\right) $ is identified with the homogeneous space $GL_n\left(
\corps\right) /Aut\left( \A\right) $. Then 
$$
\dim \vartheta \left( \A\right) =n^2-\dim Aut\left( \A\right)
$$
The orbit $\vartheta \left( \A\right) $ is provided, when $\corps=\mathbb{C}$ (a complex
field), with the structure of a differentiable manifold. In fact, $\vartheta
\left( \A\right) $ is the image through the action of the Lie group $%
GL_{n}\left( \corps\right) $ of the point $\A$, considered as a point of $%
Hom\left( V\otimes V,V\right)$. The Zariski tangent space to $Alg_{n}$ at the point $\A$
corresponds to $Z^{2}(\A ,\A)$ and the tangent space to the orbit corresponds 
to $B^{2}(\A ,\A)$.

The first approach to the study of varieties $Alg _n$ is to
establish classifications of the algebras up to isomorphisms for a fixed
dimension. Some incomplete classifications were known by mathematicians of
the last centuries : R.S. Peirce (1870),
E. Study (1890), G. Voghera (1908) and B.G. Scorza (1938). In 1974 \cite{Gabriel74}, P.
Gabriel has defined, the {\it scheme} $Alg _n$ and gave the classification,
up to isomorphisms, for $n\leq 4$  and G. Mazzola, in 1979 \cite{Mazzola79}, has studied the
case $n=5$. The number of different isomorphic classes grows up very
quickely, for example there are 19 classes in $A\lg _4$ and 59 classes in $%
A\lg _5$.

The second approach is to describe the irreducible components of a given
algebraic variety $Alg _n$. This problem has already been proposed
by Study and solved by Gabriel for $n\leq 4$ and Mazzola for $n=5.$ They
used mainly the formal deformations and degenerations.  The rigid algebras have a special interest,
an open orbit of a given algebra is dense in the irreducible component
in which it lies. Then, its Zariski closure determines an irreducible component of $Alg_{n}$, i.e. 
all algebras in this irreducible component are degenerations of the rigid algebra and there is no
algebra which degenerates to the rigid algebra. Two nonisomorphic rigid algebras correspond to different
irreducible components. So the number of rigid algebra classes gives a lower bound of the number of
irreducible components of $Alg_n$. Note that not all irreducible components 
are Zariski closure of open orbits.

 Geometrically, $\A _{0}$ is a degeneration of $\A_{1}$ means that $\A _{0}$ and $\A_{1}$ belong to the same
irreducible component in $Alg_{n}$.

The following statement gives an invariant which is stable under perturbations 
\cite{Makhlouf-Goze96} , it was used to classify the $6$-dimensional complex rigid
associative algebras and induces an algorithm to compute the irreducible components.
\begin{theorem}
Let $\A$ be a $\complex$-algebra in $Alg _n$ with a multiplication $\mu _0$ and $X_0$ be an 
idempotent of $\mu _0$. Then, for every
perturbation $\mu $ of $\mu _0$ there exists an idempotent $X$ such that $X\simeq X_0$.
\end{theorem}

This has the following consequence :
 The number of idempotents linearly independant does not decrease by perturbation.

Then, we deduce a procedure that finds the irreducible components by solving some algebraic
equations \cite{Makhlouf97}, this procedure finds all multiplication $\mu$ in $Alg _n$ with
$p$ idempotents, then the perturbations of each such $\mu$ are studied. From this it can be
decided whether $\mu$ belongs to a new irreducible component. By letting $p$ run from $n$
down to $1$, one founds all irreducible components. A computer implementation enables to do
the calculations.

Let us give the known results in low dimensions :

\begin{center}
$
\begin{array}{cccccc}
\text{dimension} & 2 & 3 & 4 & 5 & 6 \\ 
\text{irreducible\ components} & 1 & 2 & 5 & 10 & >21 \\ 
\text{rigid algebras} & 1 & 2 & 4 & 9 & 21
\end{array}
$
\end{center}

The asymptotic number of parameters in the system defining the algebraic variety $Alg_n$ is $\frac{4n^{3}}{27}+O(n^{8/3})$, see \cite{Neretin}.
Any change of basis can reduce this number by at most $n^2$ ($=dim(GL(n,\corps))$). The number of
parameters for $Alg_{2}$ is $2$ and for $Alg_{3}$ is $6$. 
 For large $n$ the
number of irreducible components $alg_{n}$ satisfies $exp(n)<alg_{n}<exp(n^{4})$, see \cite{Mazzola79}. 

Another way to study the irreducible components is the notion of compatible deformations introduced by
Gerstenhaber and Giaquinto \cite{Gerstenhaber98}. Two deformations $\A _{t}$ and  $\A _{s}$ of tha same
algebra $\A$ are compatible if they can be joined by a continuous family of algebras. When $\A$ has
finite dimension $n$, this means that  $\A _{t}$ and  $\A _{s}$ lie on a common irreducible component of
$Alg_n$. They proved the following theorem 

\begin{theorem}
Let $\A$ be a finite dimensional algebra with multiplication $\mu_{0}$ and let  $\A _{t}$ and  $\A _{s}$
 be two deformations of $\A$ with multiplications $\alpha _{t}=\mu_{0}+t F_{1}+t^{2} F_{2}+\cdot$ and 
 $\beta _{s}=\mu_{0}+s H_{1}+s^{2} H_{2}+\cdot$.
 
 If $\A _{t}$ and  $\A _{s}$ are compatible then the classes $f$ and $h$ of the cocycles $F$ and $H$
 satisfies $[f,g]_{G}=0$, i.e. $[F,H]_{G}$ must be a coboundary.
 \end{theorem}
 
 The theorem gives a necessary condition for the compatibility of the deformations. It can be used to
 show that two deformations of $\A$ lie on different irreducible components of $Alg_{n}$.


\begin{thebibliography}{99}
\bibitem{AncocheaGoze02}  Ancochea-Bermudez J. and Goze M.{\it \ On the
rigid Lie algebras, } J. Algebra 245(2002), 68-91.

\bibitem{Artin}  Artin M. {\it On the solution of analytic equations,} Inv.
Math. 5, (1968)

\bibitem{BFLS}  Bayen F., Flato M. Fronsdal C., Lichnerowicz A. and Sternheimer D. 
{\it Quantum mechanics as a deformation of classical mechanics,} 
Lett. Math. Physics 3 (1977),pp 521-530. {\it Deformation theory and quantization I and II}, Annals of
phys., 3 (1978) pp 61-110 and 111-152.

\bibitem{BMP}  Bordemann M., Makhlouf A. et Petit T. {\it D\'{e}formation par quantification et rigidit\'{e} des
alg\`{e}bres enveloppantes} Journal of Algebra, vol 285/2 (2005) pp 623-648.

\bibitem{Carles78}  Carles R. {\it Rigidit\'{e} dans la vari\'{e}t\'{e} des alg\`{e}bres}, CRASc Paris t
286 (1978), pp 1123-1226.

\bibitem{Fialowski86}  Fialowski A. {\it Deformation of Lie algebras}, Math USSR sbornik,
vol. 55 n 2 (1986), pp 467-473.

\bibitem{Fialowski88}  Fialowski A. {\it An example of formal deformations
of Lie algebras, }In Deformation theory of algebras and structures and
applications, ed. Hazewinkel and Gerstenhaber, NATO {\it \ }Adv. Sci. Inst.
Serie C 297, Kluwer Acad. Publ. (1988).

\bibitem{Fialowski90}  Fialowski A. and O'Halloran J., {\it A comparison of
deformations and orbit closure}, Comm. in Algebra 18 (12) (1990) 4121-4140.

\bibitem{Fialowski99}  Fialowski A. and Post G. , {\it versal deformations
of Lie algebra }$L_{2}$, J. Algebra 236 n 1 (2001), 93-109.

\bibitem{Fialowski00}  Fialowski A. and Fuchs D. , {\it Construction of Miniversal deformations
of Lie algebras }, J. Funct. Anal 161 n 1(1999), 76-110.

\bibitem{Fialowski02}  Fialowski A. and Schlichenmaier M. {\it Global
deformation of the Witt algebra of Krichever-Novikov type, }, Math, math/0206114, to appear in Comm.
Contemp. 

\bibitem{Gabriel74}  Gabriel P. {\it Finite representation type is open, }%
Lect. Notes in Math. 488, Springer Verlag. (1974)

\bibitem{Gerstenhaber64}  Gerstenhaber M. {\it On the deformations of rings
and algebras, }Ann. of Math 79, 84, 88 (1964, 66, 68).

\bibitem{Gerstenhaber-Schack} Gerstenhaber M. and  Schack S.D. {\it Relative Hochschild
cohomology, Rigid algebras and the Bockstein}, J. of pure and applied algebras 43
(1986), 53-74.

\bibitem{Gerstenhaber90}  Gerstenhaber M. and  Schack S.D.{\it Algebras,
bialgebras, Quantum groups and algebraic deformations }Contemporary
mathematics Vol. 134, AMS (1992), 51-92.

\bibitem{Gerstenhaber98}  Gerstenhaber M. and Giaquinto A. {\it Compatible deformation }Contemporary
mathematics Vol. 229, AMS (1998), 159-168.


\bibitem{Goze88}  Goze M. {\it Perturbations of Lie algeras structures}, In
Deformation theory of algebras and structures and applications, ed.
Hazewinkel and Gerstenhaber, NATO Adv. Sci. InsT. Serie C 297,  Kluwer Acad. Publ. (1988).

\bibitem{Goze-Makhlouf90}  Goze M. and Makhlouf A. {\it On the complex
associative algebra } Comm. in Algebra 18 (12) (1990), 4031-4046.

\bibitem{Goze-Makhlouf95}  Goze M. and Makhlouf A. {\it Elements
g\'{e}n\'{e}riques non standard en alg\`{e}bre}. Actes du colloque \`{a} la
m\'{e}moire de G. Reeb et J.L. Callot, IRMA Strasbourg, (1995).

\bibitem{Goze-Remm02}  Goze M. and E. Remm, {\it Valued deformations},
math.RA/0210475 (2002)

\bibitem{Grun_OHallo} Grunewald F. and O'Halloran J., {\it A characterization of orbit closure and
applications,} Journal of Algebra, 116 (1988), 163-175.

\bibitem{I_W} Inonu E. and Wigner E.P., {\it On the contraction of groups and their
representations,} Proc. Nat. Acad. Sci U.S., 39 (1953), 510-524.

\bibitem{Kiesler}  Keisler H.J. {\it Foundations of inifinitesimal calculus,} Prindle,
Weber and Schmidt, (1976).

\bibitem{Lutz-Goze81}  Lutz R. and Goze M. {\it Nonstandard Analysis, a
practical guide with applications,} Lect. Notes in Math. (1981).

\bibitem{lutz-makhlouf-meyer}  Lutz R, Makhlouf A et Meyer E, {\it fondement
pour un enseignement de l'Analyse en termes d'ordres de grandeur}, Brochure
APMEP n${{}^{\circ }}103$ (1996).

\bibitem{Makhlouf90}  Makhlouf A. {\it Sur les alg\`{e}bres associatives
rigides }Th\`{e}se Mulhouse (1990).

\bibitem{Makhlouf93}  Makhlouf A. {\it The irreducible components of the
nilpotent associative algebras }Revista Mathematica de la universidad
Complutence de Madrid Vol 6 n$.$1, (1993).

\bibitem{Makhlouf-Goze96}  Makhlouf A. and Goze M. {\it Classification of
rigid algebras in low dimensions}. Coll. Travaux en cours (M. Goze ed.) edition
Hermann (1996).

\bibitem{Makhlouf97}  Makhlouf A. {\it Alg\`{e}bre associative et calcul
formel }, Theoretical Computer Science 187  (1997), 123-145.

\bibitem{Massey}  Massey W.S. {\it Some higher order cohomology operations }
Symposium international de topologia algebraica, Mexico  (1958)

\bibitem{Mazzola79}  Mazzola G. {\it The algebraic and geometric
classification of associative algebras of dimension five }Manuscripta Math
27, (1979).

\bibitem{Mazzola82}  Mazzola G. {\it How to count the number of irreducible components of the
scheme of finite-dimensional algebras structures }Journal of Alg. 78, (1982) 292-302.
27, (1979).

\bibitem{Nadaud1}  Nadaud F.{\it Generalized deformations, Koszul resolutions, Moyal products}, Rev. Math. Phys.
10(5) (1998), 685-704.

\bibitem{Nelson} Nelson E. {\it Internal set theory} Annals of math (1977)

\bibitem{Neretin} Neretin Yu.A. {\it An estimate for the number of parameters defining an
$n$-dimensional algebra } Math USSR-Izv. 30 n°2 (1988), 283-294.


\bibitem{NR} Nijenhuis, A. and Richardson J.R. {\it Cohomology and deformations in
graded Lie algebras} Bull. Amer. Math. Soc 72 (1966), 1-29 

\bibitem{pincson}  Pincson G.,{\it Noncommutative deformation theory}, Lett. Math. Phys. 41 (1997), 101-117.

\bibitem{Robinson60}  Robinson A. {\it Nonstandard Analysis}, North Holland (1960)
\bibitem{schliss1} Schlessinger M. {\it Functors of Artin rings, } Trans. Amer. Math. Soc.
130 (1968),  208-222.
\bibitem{Shnider-Sternberg}  Shnider S. and Sternberg, {\it Quantum groups
from coalgebras to Drinfeld algebras}, International Press (1993).
\end{thebibliography}
\end{document}